\documentstyle[12pt]{article}\nonstopmode

\begin{document}

\newtheorem{definition}{Definition}  

\newtheorem{theorem}{Theorem}

\newtheorem{conjecture}{Conjecture}

\newtheorem{lemma}{Lemma}

\newtheorem{corollary}{Corollary}

\newtheorem{question}{Question}

\title{Locally 1-to-1 Maps and 2-to-1 Retractions}

\author{Jo Heath and Van C. Nall}
\maketitle
\begin{abstract}  This paper considers the question of which continua are 2-to-
1 retracts of continua. 

\end{abstract}
\section{Introduction.}
\hspace{.25in} A 2-to-1 retract is  a continuum that is the image of an exactly 
2-to-1 retraction defined on a continuum.

Most continua are not 2-to-1 retracts, using the word "most" as R.H. Bing did, 
because the pseudoarc is not a 2-to-1 retract; in fact, no hereditarily 
indecomposable continuum can be a 2-to-1 retract \cite{kn:h1}. Many continua 
are known not to be 2-to-1 retracts because they are not 2-to-1 images of 
continua at all. Continua in this category, excluding some that are 
hereditarily indecomposable,  include dendrites, arc-like continua, treelike 
arc continua, continua whose every subcontinuum has an endpoint,  and continua 
whose every subcontinuum has a cut point. On the other hand, if a continuum 
contains a subcontinuum that is not unicoherent then the continuum is a 2-to-1 
retract \cite{kn:nw}.  (At the end of the paper we have a glossary with 
definitions of lesser known terms.) But the fact that identifies the most  2-to-
1 retracts is that every continuum that contains a 2-to-1 retract of a 
continuum is a 2-to-1 retract \cite{kn:h3}. But note that a solenoid  shrugs 
off both criteria: a solenoid is a  2-to-1 retract, but none of its proper 
subcontinua (all arcs) are 2-to-1 retracts, and a solenoid is hereditarily 
unicoherent. 
   In Section 2 we show how to construct some 2-to-1 retracts, how to identify 
some 2-to-1 retracts and how to identify some continua that are not 2-to-1 
retracts, all using the odd fact that if some continuum $X$ maps into a 
continuum $Y$ and the map has a restriction (called a {\em 1-to-1 cover}) to an 
open proper subset of $X$ that maps 1-to-1 onto $f(X)$, then $Y$ is a 2-to-1 
retract of some other continuum.

We show, in Section 3, that maps defined on arclike continua or on hereditarily 
decomposable continua, or simple maps defined on treelike continua, have images 
that are 2-to-1 retracts provided the map is not a homeomorphism but is locally 
1-to-1 (called a {\em strictly locally 1-to-1} map). But note that the 2-to-1 
retractions themselves cannot be locally 1-to-1 at any boundary point of the 
image. We conjecture that the adjective {\em simple} (meaning $|f^{-1}(y)| \leq 
2$ for each $y $ in the image of $f$) can be removed from  the hypothesis in 
the treelike case.

In section 4 we consider decomposable continua in more detail and we show that 
if $X = A \cup B$ is a decomposable continuum and $A$ and $B$ are proper 
subcontinua of $X$, then $X$  is a 2-to-1 retract iff (1) either  $A$ or $B$ is 
a 2-to-1 retract or (2)   their intersection is not connected. This takes care 
of the decomposable case unless $A$ and $B$ can not be evaluated.

To partly justify our exclusive consideration of 2-to-1 retracts, we prove in 
the last section that if a continuum is a 2-to-1 retract of a continuum then it 
is a k-to-1 retract of a continuum for each positive integer $k$ .  
\vskip .25in
\vskip .25in
\section{Maps with 1-to-1 covers and 2-to-1 retracts}

The following theorem makes clear the connection between open covers of maps 
(see introduction or glossary for definitions) and 2-to-1 retracts, and its 
corollaries make clear its usefulness.

\vskip .25in
\begin{theorem}

The following are equivalent for the continuum $Y$:

\begin{itemize}

\item[1.] Y is a 2-to-1 retract of a continuum.

\item[2.] There is a simple map with a 1-to-1 cover from a continuum into $Y$. 

\item[3.] There is a map  with a 1-to-1 cover from a continuum into $Y$. 
\end{itemize}

\end{theorem}
Proof. Suppose $r:X \rightarrow Y$  is a 2-to-1 retraction from a continuum $X$ 
onto $ Y$.  Then $r$ is a simple map and $ U= X\setminus Y$ is an open proper 
subset of X that $ r$ maps 1-to-1 onto $ Y$.  Hence the first statement implies 
the second. And the second statement easily implies the third.
Suppose $f$ is a map from a continuum $X$ onto $Y$, and $U$ is an open proper 
subset of $X$ such that  $ f $ is 1-to-1 on $U$ and $f(U)=Y$.  Define $g:X 
\rightarrow Y \times [0, \infty)$ by $g(x) = (f(x), d(x,X \setminus U))$.  Let 
$ Y' = Y \times \{0\} $ and let $Z = Y' \cup g(X)$.  Since $Y' \cap g(X) \neq 
\emptyset$, Z is a continuum.  The 2-to-1 retraction of $ Z $ onto $ Y'$ is 
defined by $ r((y,t)) = (y,0)$. Since $Y'$ is homeomorphic to $Y$, the third 
statement implies the first.
\vskip .25in
One can construct many examples of hereditarily unicoherent 2-to-1 retracts 
using the first two corollaries to Theorem 1. For a very simple example, 
identify two points from different composants of any indecomposable continuum 
and use Corollary 2. Or use Corollary 1 and identify two disjoint subcontinua 
from different composants along a continuous map between the subcontinua. And 
Corollaries 3 and 4, rather than constructing 2-to-1 retracts, describe ways to 
decide if a given continuum is a 2-to-1 retract.
\begin{corollary}  Suppose $ X$ is a continuum, $D$ and $E$ are disjoint 
subcontinua of $X$, and $h$ is a map from D into E. Then $Y = X/\{\{x,h^{-
1}(x)\} | x \in E\}$ is a 2-to-1 retract of a continuum.

\end{corollary}
Proof. Let $U$ be $X \setminus D$; then the quotient map $p:X \rightarrow Y$ 
maps $U$ 1-to-1 onto $Y$.
\begin{corollary}  Suppose $ X$ is a continuum and $p$ and $q$ are two points 
of $X$. Then $X/\{p,q\}$ is a 2-to-1 retract of a continuum.

\end{corollary}

\begin{corollary} Suppose $Y$ is a continuum and $K$ is a local cut continuum 
that is not a cut continuum, i.e. $ Y \backslash K $ is connected but there is 
an open set $U$ containing $K$ such that $ U \backslash K = A \cup B$, two 
nonempty separated sets, and $K$ contains both a limit point of $A$ and a limit 
point of $B$. Then $Y$ is a 2-to-1 retract of a continuum.

\end{corollary}
Proof. Construct a continuum $X$ by adding to the connected set $Y \backslash 
K$ two disjoint copies of $K$, say $K^{1}$ and $K^{2}$, with $K^{1}$ attached 
to $A$ in the same way that $K$ was attached  to $A$ and with $K^{2}$ attached 
to $B$ in the same  way that $K$ was attached  to $B$. Let $h$ be the 
homeomorphism from $K^{1}$  to $ K^{2}$ such that for each point $t$ in $K$, 
$h$ takes the copy of $t$ in $K^{1}$ to the copy of $t$ in $K^{2}$. Then, by 
Corollary 1, $Y = X/\{\{x,h^{-1}(x)\} | x \in  K^{2}\}$ is a 2-to-1 retract of 
a continuum.

\begin{corollary} If the continuum $Y$ has a local cut point that is not a cut 
point, then $Y$ is a 2-to-1 retract.

\end{corollary}
\vskip .25in
\section{Strictly locally 1-to-1 maps and 2-to-1 retracts}
The next series of results are intended to demonstrate that the strictly 
locally 1-to-1 image of a continuum is frequently a 2-to-1 retract because it 
has a 1-to-1 cover.  Later we have two examples that demonstrate the sort of 
complexity that a continuum might have in order for it to have a strictly 
locally 1-to-1 image that is not a 2-to-1 retract.

\begin{lemma}

If $f$ is strictly locally 1-to-1 map from the continuum $X$ into the continuum 
$Y$, and $f$ is 1-to-1 on the closed subset $A$ of $X$ and 1-to-1 on $X 
\backslash A$, then $Y$ is a 2-to-1 retract.

\end{lemma}
Proof. Let $X_o = \{ x \in X \mid f^{-1}(f(x)) \neq \{x \} \}$.  Since $f$ is 
locally 1-to-1 $X_o$ is closed.   Let $U = X \backslash (X_o \cap A)$, and  $ 
U$ is clearly an open set.  Since $f$ is 1-to-1 on $A$, for each $x$ in $X_o 
\cap A$ there is an $\hat{x}$ in $U$ such that $f(x) = f(\hat{x})$.  Therefore 
$f(U) = f(X)$.  There cannot be three elements of $X$ with the same image under 
$f$ since $f$ is 1-to-1 on $A$ and on $X \backslash A$, and if there are two 
elements of $X$ with the same image under $f$, then one of them is in $X_o \cap 
A$.   Therefore $f$ is 1-to-1 on $U$. It follows now from Theorem 1 that $Y$ is 
a 2-to-1 retract.

\vskip .25in

\begin{lemma}

 If $f$ maps the compactum $X$ onto $Y$ so that (1) $f$ is strictly locally 1-
to-1, (2) $f$ is 1-to-1 on each proper subcontinuum of $X$, and (3) there is at 
least one 1-to-1 point $p$ in $X$  (meaning that no other point in $X$ maps to 
$f(p)$), then $f$ has a 1-to-1 cover.

 \end{lemma}
Proof. First, for each set $K$ in $ X $, define $ \hat{K} $ to be the points in 
$ X  \setminus K $ that map the same under $ f $ as some point in $ K $. Note 
that if $ K $ is closed, then so is $ \hat{K} $. Since the set of 1-to-1 points 
in $ X $ is open, there is an open set $ U $ containing $ p $ that is contained 
in the set of 1-to-1 points. The components of $ X \setminus U $ are components 
of a compactum, so if $ C $ is such a component and $ \epsilon > 0 $, then 
there is an open and closed set $ V(C) $ in the $ \epsilon $-neighborhood of $ 
C $ that contains $ C $; further, since $ f $ is locally 1-to-1 and 1-to-1 on 
each subcontinuum, we may assume that $ f $ is 1-to-1 on $ V(C) $. $V(C) $ is 
open and closed in $ X \setminus U $. Let $ V_{1}, V_{2}, ..., V_{n} $ be a 
finite cover of $ X \setminus U $ consisting of these $ V(C) $ sets. Now, let 
$$ W = U  \cup V_{1} \cup (V_{2} \setminus \hat{V_{1}}) \cup ... \cup (V_{n} 
\setminus (\hat{V_{1}} \cup \hat{V_{2}} \cup ... \cup \hat{V}_{(n-1)})).$$
Each $ \hat{V_{i}} $ is closed, so the parenthetical sets are each open (in $ 
X\setminus U $). Since $ W \setminus U $ is open in $ X \setminus U $, $ W $ is 
open in $ X $. And $ f $ is 1-to-1 on $ W $ and maps $ W $ onto $Y$, so $f|W$ 
is a 1-to-1 cover of $f$.

\begin{corollary} If $f$ is a strictly locally 1-to-1 map from a continuum $X$ 
into a continuum $Y$, $f$ is 1-to-1 on each proper subcontinuum of $X$, and 
there is at least one 1-to-1 point for $f$, then $Y$ is a 2-to-1 retract.

\end{corollary}
\begin{lemma}

If $f$ is a strictly locally 1-to-1 map from  a decomposable continuum $X$ into 
a continuum $Y$, and $f$ is 1-to-1 on each proper subcontinuum of $X$, then $Y$ 
is a 2-to-1 retract.

\end{lemma}

Proof. Since $X$ is decomposable, $X$ is the union of two proper subcontinua, 
$A$ and $B$; and since  $f$ is 1-to-1 on each proper subcontinuum of $X$, $f$ 
restricted to each of $A$ and $B$ is 1-to-1. Thus every point of $A \cap B$ is 
a 1-to-1 point and the hypothesis of Corollary 5 is satisfied. Hence $Y$ is a 2-
to-1 retract of a continuum. 

\vskip .25in
\begin{theorem} 

If $f$ is a strictly locally 1-to-1 map from a hereditarily decomposable 
continuum $X$ into a continuum $Y$, then $Y$ is a 2-to-1 retract of a 
continuum.

\end{theorem}
Proof. If $X'$ is minimal with respect to being a subcontinuum of $X$ on which 
$f$ is not 1-to-1 then the conditions of the previous lemma are satisfied by 
the restriction of $f$ to $X'$. So $f(X')$ is a 2-to-1 retract of a continuum 
by Lemma 3 and, since every continuum that contains a 2-to-1 retract is itself 
a 2-to-1 retract, $f(X)$ is also a 2-to-1 retract of a continuum.

\vskip .25in
\begin{theorem}

The image of a strictly locally 1-to-1 map defined on an arc-like continuum is 
a 2-to-1 retract of a continuum.

\end{theorem}
Proof. Assume $ X$ is an arc-like continuum, and $ f$ is a strictly locally 1-
to-1  map with domain $X$.  Since $ f$ is locally 1-to-1, there is a positive 
number $\epsilon$ such that if $ f(x) = f(y)$ and $x \neq y$, then $d(x,y) > 
\epsilon$.  Let g be an $\epsilon$-map onto $[0,1]$, and let $A= \{ x \in X 
\mid \exists x' \neq x  \ni f(x')=f(x)$ and $g(x') < g(x)\}$.  It is easy to 
verify that $A$ is closed, and that if $U= X \backslash A$, then $f$ is 1-to-1 
on $U$ and $f(U)=f(X)$. Hence $f$ has a 1-to-1 cover and its image must be a 2-
to-1 retract of a continuum.
\begin{lemma}

If the continuum $X$ is the union of two continua $A$ and $B$ and every 
strictly locally 1-to-1 image of $A$ and every strictly locally 1-to-1 image of 
$B$ is a 2-to-1 retract of a continuum, then every strictly locally 1-to-1 
image of $X$ is a 2-to-1 retract of a continuum.

\end{lemma}
Proof. A strictly locally 1-to-1 map with domain $X$ is either strictly locally 
1-to-1 on $A$, strictly locally 1-to-1 on $B$, or 1-to-1 on $A$ and on  
$X\backslash A$.  In each case $f(X)$ is a 2-to-1 retract of a continuum; in 
the latter case by Lemma 1 and in the first two cases because $f(X)$ contains a 
2-to-1 retract of a continuum.

\vskip .25in

\begin{theorem}

If $f$ is a strictly locally 1-to-1 map defined on a  continuum $X$ that is a 
finite union of continua which are either arc-like or hereditarily 
decomposable, then $f(X)$ is a 2-to-1 retract of a continuum.

\end{theorem}
Proof. Theorem 4 follows from Lemma 4, Theorem 2 and Theorem 3.
\vskip .25in

We would like to be able to replace arc-like with tree-like in  Theorem 3.  In 
Theorem 5 we come close, but there is an added assumption that the map is 
simple.  We conjecture that this assumption is not necessary.

\vskip .25in
\begin{lemma}

No tree-like continuum admits a non-trivial $k$-fold covering map.

\end{lemma}
Proof.  Every $k$-fold covering map is open and therefore, by a theorem of 
G.T.  Whyburn \cite[Theorem 7.5, p. 148]{why}, confluent.  McLean \cite{mclean} 
has shown that the confluent image of a tree-like continuum is itself a tree-
like continuum and Ma\'{c}kowiak \cite{mack} has shown that a local 
homeomorphism onto a tree-like continuum is a homeomorphism.  Hence, any 
covering map defined on a tree-like continuum must be the trivial 1-to-1 
covering map.
\vskip .25in
\begin{theorem}

The image of a simple strictly locally 1-to-1 map defined on  a treelike 
continuum is a 2-to-1 retract of a continuum.

\end{theorem}
Proof. Suppose we have a simple strictly locally 1-to-1 map defined on  a 
treelike continuum; then there is a restriction, say $f$, to a tree-like 
subcontinuum $ X $ of the domain such that  $f$ is strictly locally 1-to-1 and 
is 1-to-1 on each proper subcontinuum of $X$. So $f$ cannot be a covering map 
by the previous lemma. Hence, since it is locally 1-to-1 it cannot be exactly 2-
to-1; and so, since $ f $ is simple, there is a point in $ X $ at which $f$ is 
1-to-1. Thus, by Corollary 5, $f(X)$, and thus the original image space, is a 2-
to-1 retract of a continuum. 
\vskip .25in
\begin{question} Is the hypothesis that the map be simple necessary in Theorem 
5?

\end{question}

\vskip .25in
\begin{corollary}

If $f$ is a simple strictly locally 1-to-1 map defined on a continuum $X$ that 
is a finite union of continua that are either tree-like or hereditarily 
decomposable, then $f(X)$ is a 2-to-1 retract.

\end{corollary}
\vskip .25in
To find an example of a continuum that has a strictly locally 1-to-1 image that 
is not a 2-to-1 retract it is natural to think of a continuum that is 
hereditarily indecomposable  with a locally 1-to-1 image that is also 
hereditarily indecomposable.  That makes a 2-fold cover from the  pseudo-circle 
onto itself a natural choice.  Note that the pseudo-circle is an example of a 
continuum that is a 2-to-1 image of a continuum but is not a 2-to-1 retract of 
a continuum. In the second example the domain and range are decomposable, but 
just barely so.
\vskip .25 in
\noindent {\bf Example 1.} A pseudo-circle is a hereditarily indecomposable, 
circularly chainable, separating plane continuum.  It was shown in 
\cite[Example 1]{kn:h2} that there is a 2-fold cover, and therefore a strictly 
locally 1-to-1 map, from the pseudo-circle onto itself , and in \cite[Theorem 
5]{kn:h1} that no hereditarily indecomposable continuum is a 2-to-1 retract of 
a continuum.  The 2-fold cover is a simple strictly locally 1-to-1 map but 
every restriction of the 2-fold cover to a proper subcontinuum of the pseudo-
circle is a homeomorphism.

\vskip .25in

\noindent {\bf Example 2.} The continuum $X$ is the union of two pseudo-
circles, $P_{1}$ and $P_{2}$, joined at two points, and its image $Y$ is the 
union of two pseudo-circles, $Q_{1}$ and $Q_{2}$, joined at one point. As 
mentioned in Example 1, there are 2-fold covers, $g_{1}$ and $g_{2}$, from 
$P_{1}$ onto $Q_{1}$, and from $P_{2}$ onto $Q_{2}$, respectively. Suppose $a$ 
and $b$ are points in $P_{1}$ such that $g_{1}(a) = g_{1}(b) $, and $c$ and $d$ 
are points in $P_{2}$ such that $g_{2}(c) = g_{2}(d) $.  To form $X$, attach 
$a$ in $P_{1}$ to  $c$ in $P_{2}$, and  attach  $b$ in $P_{1}$ to  $d$ in 
$P_{2}$. To form $Y$, attach $g_{1}(a)$ in $Q_{1} = g_{1}(P_{1})$ to  
$g_{2}(c)$ in $Q_{2} = g_{2}(P_{2})$. Then the map $g_{1} \cup g_{2}$ is a 
simple, strictly locally 1-to-1 map from $X$ onto $Y$.  Since the pseudo-circle 
is not a 2-to-1 retract, $Y$ is not a 2-to-1 retract by Theorem 6 which is 
proven below.

\vskip .25in
\section{When are decomposable continua 2-to-1 retracts?}
Suppose $X = A \cup B$ is a decomposable continuum, and $A$ and $B$ are proper 
subcontinua. When is $X$ a 2-to-1 retract? If $A \cap B$ is not connected then 
$X$ is not unicoherent and we know from \cite{kn:nw} that $X$ is a 2-to-1 
retract. If either $A $ or $B$ is a 2-to-1 retract, then we know from 
\cite{kn:h3} that $X$ is a 2-to-1 retract. But, are these conditions necessary? 
Yes. We show in Theorem 6 that if $A$ and $B$ both fail to be 2-to-1 retracts 
and if their intersection is connected, then $X$ cannot be a 2-to-1 retract.
\begin{lemma} If $X$ is a 2-to-1 retract, and $K$ is a subcontinuum  of $X$, 
then there is a component $C$ of $X \backslash K$ such that $C \cup K $ is a 2-
to-1 retract.

\end{lemma}
Proof. Let $r: Z \rightarrow X$ be a 2-to-1 retraction from a continuum $Z$ 
onto $X$. If $r^{-1}(K)$ is connected, then $K$ is a 2-to-1 retract, so the 
conclusion is true for any component of $X \backslash K$.  So, assume that $r^{-
1}(K)$ is not connected. Then $r^{-1}(K)$  is contained in $D \cup E$, where 
$D$ and $E$ are disjoint open sets intersecting $r^{-1}(K)$. Without loss of 
generality, we will assume that $K$ is in $D$.  Let $K'$ be a copy of $K$ 
disjoint from $Z$. For each point $x$ in $r^{-1}(K) \cap D$, identify $x$ with 
$r(x)$, and for each point $x$ in $r^{-1}(K)  \cap E$, identify $x$ with the 
copy of $r(x)$ in $K'$. Call this new continuum $Z'$. We then have a 2-to-1 
retraction $R$ from $Z'$ onto $X$ for which $R^{-1}(K)$  has exactly two 
components, $K$ and $K'$. 
There is a component $C'$ of $Z' \backslash (X \cup K')$ whose closure 
intersects both $X$ and $K'$ since $Z'$ is connected. Let $C$ denote the 
component of $X \backslash K$ that contains $R(C')$. Some point $x$ of $X$ is 
the limit of a sequence $S$ of points of $C'$ and $x$ must also be the limit of 
the sequence $R(S)$. Hence $T = K \cup C \cup C' \cup K'$ is connected. Every 
component of $Z' \backslash (X \cup K')$ either maps into $C$ or its image 
misses $C$, and the closure of each component of $Z' \backslash (X \cup K')$ 
intersects either $K'$ or $X$. Suppose such a component $V$ maps into $C$ . If 
its closure intersects $K'$, then  $V \cup K' $ is connected and if its closure 
intersects $X$ then its closure intersects $C$ by the same argument that the 
closure of $C'$ intersects $C$, so $V \cup C $ is connected. Hence, all of the 
components of $Z' \backslash (X \cup K') $ that map into $C$ can be added to 
$T$, getting a connected set that is equal to $R^{-1}(K \cup C)$. Thus $K \cup 
C$ is a 2-to-1 retract.
\begin{theorem} Suppose $X = A \cup B$ is a decomposable continuum and each of 
$A$ and $B$ is a proper subcontinuum. Then $X$ is a 2-to-1 retract iff one of 
the following is true:
\begin{itemize}
\item $A$ is a 2-to-1 retract, or
\item $B$ is a 2-to-1 retract, or
\item $ A \cup B$ is not connected.
\end{itemize}
\end{theorem}
Proof.  The sufficiency of each of the three conditions is discussed at the 
beginning of Section 4. For the converse, assume $X = A \cup B$ is a 2-to-1 
retract of a continuum, and that $K = A \cap B $ is connected. Then, by Lemma 
6, there is a component $C$ of $X \backslash K $ such that $C \cup K$ is a 2-to-
1 retract. But $C $ must either be a subset of $A$ or of $B$. If $C \subset A$, 
then $C \cup K$ is a 2-to-1 retract in $A$ which implies that $A$ itself is a 2-
to-1 retract. Thus, one of the three conditions has to hold.
\vskip .25in
\section{2-to-1 retract implies k-to-1 retract.}
Information we have on which continua are 2-to-1 retracts helps with the study 
of which continua are k-to-1 retracts, for other positive integers k, by way of 
the corollary below.
\begin{theorem} 

Suppose $Y$ is a $k$-to-1 retract of a continuum. Then, for each positive 
integer $n$, $Y$ is a $(1 + (k-1) \times n )$-to-1 retract of a continuum.

\end{theorem}
Proof. Suppose $X$ is a continuum and $r:X \rightarrow Y$ is a k-to-1 
retraction onto $Y$.     Let $n$ be a positive integer.  Define the map $g_{i} 
: X \rightarrow X \times \prod_{j=1}^{n}{[0,\infty)} $ for $1 \leq i < n$ by 
letting $g_{i}(x)$ be the point in  $X \times \prod_{j=1}^{n}{[0,\infty)}$ with 
first coordinate $ x$, with $ i+1$ coordinate $d(x,Y)$, and  with all other 
coordinates zero.  Let $Y' = Y \times \prod_{j=1}^{n}{\{0\}}$.  Let $ Z = Y' 
\cup  (\bigcup _{i=1}^{n} {g_{i} (X)})$.  Since $g_{i}(X)$ intersects $Y'$ for 
each $i$, $Z$ is a continuum.  The $(1 + ((k-1) \times n))$-to-one retraction 
$r^* : Z \rightarrow Y'$ is defined by $r^*((x,t_1,t_2, . . . , t_n)) = 
(r(x),0,0,. . . , 0)$. The conclusion of the theorem follows because $Y$ is 
homeomorphic to $Y'$.
\vskip .25 in
\begin{corollary}

If a continuum $Y$ is a 2-to-1 retract of a continuum, then $Y$ is a $k$-to-1 
retract of a continuum, for each $k > 2$.

\end{corollary}

\section{Glossary}
\begin{itemize}
\item {\bf Arclike.} A continuum is {\em arclike} if for each $\epsilon > 0$, 
there is an $\epsilon$-map from the continuum onto an arc.

\item {\bf Confluent Map}. A map is {\em confluent} if each component of the 
inverse of any continuum $C$ in the image is mapped onto $C$.

\item {\bf Continuum.}  A topological space is a {\em continuum} if it is 
connected, compact, and metric.
\item {\bf Covering Map.} A map defined on a continuum is a {\em covering map} 
if it is k-to-1 for some positive integer k, open, and locally 1-to-1.
\item {\bf Indecomposable Continuum.} A continuum is {\em indecomposable} if it 
is not the union of two proper subcontinua.
\item {\bf Local Cut Continuum and Local Cut Point}  A subcontinuum $K$ of a 
continuum $Y$ is a local cut continuum if $ Y \backslash K $ is connected but 
there is an open set $U$ containing $K$ such that $ U \backslash K = A \cup B$, 
two nonempty separated sets, and $K$ contains both a limit point of $A$ and a 
limit point of $B$. If $K$ consists of a single point then that point is called 
a local cut point. 
\item {\bf Map.} A {\em map} is a continuous function.
\item {\bf 1-to-1 cover} A {\em 1-to-1 cover} of a map $f$ with domain $X$ is a 
restriction of $f$ to an open proper subset $U$ of $X$  such that $f$ is 1-to-1 
on $U$ and $f(U) = f(X)$. 
\item {\bf Simple Map} A map is {\em simple} if the cardinality of each point 
inverse is either one or two.
\item {\bf Strictly locally 1-to-1.} A {\em strictly locally 1-to-1 map} is a 
map which is locally 1-to-1 but not a homeomorphism.
\item {\bf Treelike.} A continuum is {\em treelike} is for each $\epsilon > 0$,

there is an $\epsilon$-map from the continuum onto a tree (an acyclic graph).
\item {\bf 2-to-1.} A function is {\em 2-to-1} if the preimage of each point in 
the image has exactly two points.
\item {\bf 2-to-1 retract} A continuum $Y$ is a {\em 2-to-1 retract} if there 
is a continuum $X$ and a retraction $r$ from $X$ onto a subcontinuum of $X$ 
that is homeomorphic to $Y$.
\item {\bf Unicoherent Continuum} A continuum is {\em unicoherent} if it is not 
the union of two subcontinua whose intersection is not connected.
\end{itemize}


\begin{thebibliography}{V1} 


\bibitem{kn:h1} Jo Heath, {\em A non-treelike continuum that is not the 2-to-1 
image of any continuum }, Proceedings of the AMS {\bf 124} (1996) 3571-3378.

\bibitem{kn:h2} Jo Heath, {\em Weakly confluent, 2-to-1 maps on hereditarily 
indecomposable continua}. Proceedings of the AMS {\bf 117} (1993) 569-573.

\bibitem{kn:h3} Jo Heath, {\em Exactly k-to-1 functions: from pathological 
functions with finitely many discontinuities to well-behaved covering maps} 
Continua with the Houston problem book, Marcel Dekker, Inc. New York, 89 - 103.

\bibitem{kn:nw}  S. B. Nadler, Jr. and L. E. Ward, Jr., {\em Concerning exactly 
(n,1) images of continua}  Proceedings of the AMS {\bf 87} (1983), 351-354.

\bibitem{mack} T. Ma\'{c}kowiak, {\em Local homeomorphisms onto tree-like 
continua}, Colloq. Math., {\bf 38} (1977), 63-68.

\bibitem{mclean} R. Bruce McLean, {\em Confluent images of tree-like curves are 
tree-like}, Duke J. {\bf 39} (1972), 465-473.

\bibitem{why} G.  T.  Whyburn, {\em Analytic topology}, Amer.  Math.  Soc. 
Colloq.  Publ., {\bf 28}, Amer.  Math.  Soc., Providence, R.I., 1963

\end{thebibliography}
\end{document}